\documentclass[preprint,11pt]{elsarticle}

\usepackage{amssymb}
\usepackage{amsmath}
\usepackage{amsfonts}
\usepackage{amssymb}
\usepackage{indentfirst,latexsym,bm}
\usepackage{amsthm}
\usepackage{color,xcolor}
\usepackage[all]{xy}
\input{mathrsfs.sty}
\newtheorem{theorem}{Theorem}[section]
\newtheorem{lemma}[theorem]{Lemma}
\newtheorem{corollary}[theorem]{Corollary}

\theoremstyle{definition}

\theoremstyle{definition}

\theoremstyle{remark}
\newtheorem{remark}{Remark}[section]
\theoremstyle{question}

\theoremstyle{problem}
\newtheorem{problem}{Problem}[section]

\numberwithin{equation}{section}

\journal{XXX}

\begin{document}

\begin{frontmatter}






\title{Convex inequalities in Hilbert $C^*$-modules}
\author[shnu]{Kangjian Wu}
\ead{wukjcool@163.com}
\author[shnu]{Jia Li}
\ead{1184048619@qq.com}
\author[shnu]{Qingxiang Xu}
\ead{qingxiang\_xu@126.com}
\address[shnu]{Department of Mathematics, Shanghai Normal University, Shanghai 200234, PR China}
\begin{abstract}The H$\ddot{{\rm o}}$lder-McCarty inequalities are originally derived in the Hilbert space case and have been generalized via a convex inequality. The main purpose of this paper is to extend this convex inequality to the Hilbert $C^*$-module case, and meanwhile to
make some investigations on the H$\ddot{{\rm o}}$lder-McCarty inequalities in the Hilbert $C^*$-module case.



\end{abstract}


\end{frontmatter}



\section{Introduction and the preliminaries}\label{sec:pre}

Given a Hilbert space $H$, let $\mathbb{B}(H)$ be the set of all bounded linear operators on $H$.
The identity and the positive elements in the $C^*$-algebra $\mathbb{B}(H)$ are denoted by $I_H$ and $\mathbb{B}(H)_+$, respectively. Let $T\in\mathbb{B}(H)_+$ and $x\in H\setminus\{0\}$ be arbitrary. The H$\ddot{{\rm o}}$lder-McCarty inequalities (see \cite[Lemma~1.3.3]{BDMP} and \cite[Lemma~2.1]{McCarthy}) are known as
\begin{align*}&\langle Tx,x\rangle^r\leq\|x\|^{2(r-1)}\langle T^rx,x\rangle, \quad r\ge 1,\\
&\langle Tx,x\rangle^r\geq\|x\|^{2(r-1)}\langle T^rx,x\rangle,\quad 0<r\le 1,
\end{align*}
which can be stated equivalently as
\begin{align*}&\langle Tx,x\rangle^r\leq \langle T^rx,x\rangle, \quad r\ge 1,\\
&\langle Tx,x\rangle^r\geq \langle T^rx,x\rangle,\quad 0<r\le 1
\end{align*}
for every $T\in\mathbb{B}(H)_+$ and $x\in H$ with $\|x\|=1$. These inequalities are particularly useful in dealing with the numerical radius \cite{BDMP,KR,MRO,Rashid,RA,Yang-Li}. Note that when $T\in\mathbb{B}(H)_+$ and $r\ge 1$, the function $f(u)=u^r$ is  continuous and convex on the closed interval $[0,\|T\|]$. The same is true if $0<r\le 1$ and the function is replaced by $f(u)=-u^r$. So, some attention has been paid in the literature to study convex inequalities. As a generalization of the H$\ddot{{\rm o}}$lder-McCarty inequalities, a convex inequality was provided by Mond and Pe$\breve{{\rm c}}$ari$\acute{{\rm c}}$ \cite{MP} as follows.

\begin{lemma}\label{lem:convex equal}{\rm \cite[Theorem~1]{MP}}
Let $\mathbb{R}$ denote the set of the real numbers. If $m,M\in\mathbb{R}$ and $T\in\mathbb{B}(H)$ is self-adjoint such that $mI_H\le T\le MI_H$, then for every continuous convex function $f: [m,M]\to \mathbb{R}$,
we have
\begin{equation*}f\big(\langle Tx, x\rangle\big)\leq \big\langle f(T)x,x \big\rangle,\end{equation*}
where $x$ is an arbitrary element of $H$ satisfying $\|x\|=1$.
\end{lemma}

Mond and Pe$\breve{{\rm c}}$ari$\acute{{\rm c}}$ used the following elementary result in \cite{MP} without a proof, which plays a crucial role in their proof Lemma~\ref{lem:convex equal}.

\begin{lemma}\label{lem:an elementray observation}Let $f: [a,b]\to\mathbb{R}$ be a continuous convex function\footnote{Every convex function on $[a,b]$ is continuous on $(a,b)$. It may however fail to be continuous at the endpoints $a$ and $b$.}. Then for every $x_0\in [a,b]$ and  $\varepsilon>0$, there exists
a function $l(x)$ given by
$$l(x)=cx+d \quad (x\in \mathbb{R})$$ for some constants $c$ and $d$ such that the following conditions are satisfied:
\begin{itemize}
\item[{\rm (i)}] $f(x)\ge l(x)$ for all $x\in [a,b]$;
\item[{\rm (ii)}] $f(x_0)<l(x_0)+\varepsilon$.
\end{itemize}
\end{lemma}

For the sake of completeness, in this paper we will give a detailed proof of Lemma~\ref{lem:an elementray observation}. The main purpose of this paper is to extend Lemma~\ref{lem:convex equal} to the Hilbert $C^*$-module case, and meanwhile to
make some investigations on the H$\ddot{{\rm o}}$lder-McCarty inequalities in the Hilbert $C^*$-module case.

Let us recall some basic knowledge about (right) Hilbert $C^*$-modules. Hilbert $C^*$-modules are the natural generalizations of Hilbert spaces by allowing the inner product to take values in a $C^*$-algebra rather than in the complex field $\mathbb{C}$ \cite{Lance,MT,Paschke}. Given Hilbert $C^*$-modules $H$ and $K$ over a $C^*$-algebra $\mathfrak{A}$, let $\mathcal{L}(H,K)$ denote the set of all adjointable operators from $H$ to $K$, with the abbreviation $\mathcal{L}(H)$ whenever $H=K$. Let $\mathcal{L}(H)_{\mbox{sa}}$ (resp.\,$\mathcal{L}(H)_+$) denote the set of all self-adjoint (resp.\,positive) elements in the $C^*$-algebra $\mathcal{L}(H)$. The unit of $\mathcal{L}(H)$, i.e., the identity operator on $H$, is denoted by $I_H$. In the case that $A\in\mathcal{L}(H)_+$, the notation $A\ge 0$ is also used to indicate that $A$ is a positive operator on $H$.

It is known that every adjointable operator $A\in\mathcal{L}(H,K)$ is a bounded linear operator, which is also $\mathfrak{A}$-linear in the sense that
\begin{equation}\label{equ:keep module operator}A(xa)=A(x)a,\quad\forall\,x\in H,  a\in\mathfrak{A}.
\end{equation}
Let $\mathcal{R}(A)$ and $A^*$ denote the range and the adjoint operator of $A$, respectively.

Unless otherwise specified, throughout the rest of this paper $\mathbb{N}$ is the set consisting of all positive integers, $M_n(\mathbb{C})$ stands for the set of all $n\times n$ complex matrices,  $\mathfrak{A}$ is a $C^*$-algebra, $\mathcal{S}(\mathfrak{A})$ is the set of all states on $\mathfrak{A}$, and $E$ is a  Hilbert $\mathfrak{A}$-module. For every $T\in\mathcal{L}(E)$, we take the convention that $T^0=I_E$.

A brief description of the paper is as follows. A detailed proof of Lemma~\ref{lem:an elementray observation} is provided in Section~\ref{sec:proof of elementray lemma}. Section~\ref{sec:generalization} aims to extend Lemma~\ref{lem:convex equal} to the Hilbert $C^*$-module case. Some investigations on the H$\ddot{{\rm o}}$lder-McCarty inequalities are dealt with in Section~\ref{sec:investigation}.
Specifically, the Hilbert $C^*$-module version of the H$\ddot{{\rm o}}$lder-McCarty inequalities is provided in Theorem~\ref{thm:generalized McCarthy inequality}. In addition, two problems (Problems~\ref{prob01} and \ref{prob02}) are raised, and two partial solutions are provided in this section; see Theorems~\ref{thm:commutative C-star alg case} and  \ref{thm:negative answer} for the details.

\section{The proof of Lemma~\ref{lem:an elementray observation}}\label{sec:proof of elementray lemma}
Several cases are taken into consideration.

\textbf{Case 1}: $x_0\in (a,b)$. In this case, both of the left-hand derivative $f_{-}^{'}(x_0)$ and the right-hand derivative $f_{+}^{'}(x_0)$ exist such that $f_{-}^{'}(x_0)\leq f_{+}^{'}(x_0)$.

\textbf{Subcase 1.1}: If $x \in (x_0, b]$, then $\frac{f(x)-f(x_0)}{x-x_0}\geq f_{+}^{'}(x_0)$, so
\begin{align*}
 f(x) & \geq f_{+}^{'}(x_0)(x-x_0)+f(x_0)\\
      & \geq f_{-}^{'}(x_0)(x-x_0)+f(x_0).
\end{align*}

\textbf{Subcase 1.2}: If $x \in [a,x_0)$. then $\frac{f(x)-f(x_0)}{x-x_0}\leq f_{-}^{'}(x_0)$. Thus,
\begin{align*}
 f(x) \geq f_{-}^{'}(x_0)(x-x_0)+f(x_0).
\end{align*}

Combining Subcases 1.1 and 1.2 yields
\begin{align*}
 f(x) \geq f_{-}^{'}(x_0)(x-x_0)+f(x_0), \quad\forall x \in [a,b].
\end{align*}
Let $c=f_{-}^{'}(x_0), d=-x_0f_{-}^{'}(x_0)+f(x_0)$, and $l(x)=cx+d$. Then $l(x)$ satisfies conditions (i) and (ii) of this lemma.

\textbf{Case 2}: $x_0=a$. In this case,  $f_{+}^{'}(a) \in\mathbb{R}$ or $f_{+}^{'}(a)=-\infty$.

\textbf{Subcase 2.1}: $f_{+}^{'}(a)\in\mathbb{R}$. For every $x\in (a,b]$, we have $\frac{f(x)-f(a)}{x-a}\geq f_{+}^{'}(a)$, which gives
$f(x) \geq f_{+}'(a)(x-a)+f(a)$. Therefore, $l(x)$ satisfies conditions (i) and (ii) of this lemma, in which
$l(x)=cx+d$ with $c=f_{+}'(a)$ and $d=-af_{+}'(a)+f(a)$.

\textbf{Subcase 2.2}: $f_{+}^{'}(a) =-\infty$. Choose $\delta_1>0$ such that $a+\delta_1 < b$ and
$$\frac{f(x)-f(a)}{x-a}<0,\quad \forall x\in (a,a+\delta_1].$$ So $f(x)-f(a )<0$ for every $x\in (a,a+\delta_1]$, which indicates that $f(a)$ is the maximum value of $f$ on $[a,a+\delta_1]$. For any $\varepsilon>0$, since $f$ is continuous from the right  at $x=a$, there exists $\delta \in (0, \delta_1)$ such that
\begin{align*}
    |f(x)-f(a)|< \frac{\varepsilon}{2},\quad \forall x \in [a, a+\delta].
\end{align*}
As $f$ is continuous on the finite closed interval $[a,a+\delta]$, it takes the  minimum value at some point $x_1$ in $[a, a+\delta]$. If $x_1=a$, then
$f$ is constant on $[a, a+\delta]$ (since $f(a)$ is the maximum value of $f$ on this interval). This contradicts the assumption that $f_{+}^{'}(a) =-\infty$.
Therefore, $x_1>a$ and thus the notation $f'_{-}(x_1)$ is meaningful such that
\begin{align*}
    f'_{-}(x_1)=\lim_{x\rightarrow{x_1^-}}\frac{f(x)-f(x_1)}{x-x_1} \leq 0.
\end{align*}

Now let $l(x)=f'_{-}(x_1)(x-x_1)+f(x_1)=cx+d$, in which $c=f'_{-}(x_1)$ and $d=-x_1f'_{-}(x_1)+f(x_1)$. Then $l(x_1)=f(x_1)$, and as shown in Case 1 we have $f(x) \geq l(x)$ for all $x \in [a,b]$.
Furthermore,
\begin{align*}f(x_0)&=f(a)=f(a)-f(x_1)+f(x_1)=|f(x_1)-f(a)|+l(x_1)\\
&<\frac{\varepsilon}{2}+l (x_1) =\frac{\varepsilon}{2} + l(a) + l(x_1)-l(a)=\frac{\varepsilon}{2} + l(x_0)+c(x_1-a)\\
&\le \frac{\varepsilon}{2} + l(x_0) <\varepsilon + l(x_0).
\end{align*}
This also shows the validity of (i) and (ii) in this lemma.

\textbf{Case 3}: $x_0=b$. In this case, either $f'_{-}(b)\in\mathbb{R}$ or $f'_{-}(b)=+\infty$.

\textbf{Subcase 3.1}: $f'_{-}(b)\in\mathbb{R}$. For every $x \in [a,b)$, we have
$\frac{f(x)-f(b)}{x-b} \leq f'_{-}(b)$, which yields
\begin{align*}
     f(x)\geq f'_{-}(b)(x-b)+f(b).
\end{align*}
Hence, the conclusion holds if we put $l(x) = cx+d$ with  $c=f'_{-}(b)$ and $d=-bf'_{-}(b)+f(b)$.

\textbf{Subcase 3.2}: If $f'_{-}(b) = +\infty$. Choose $\delta_1>0$ such that $a <b-\delta_1$ and $\frac{f(x)-f(b)}{x-b}>0$ for every $x \in [b-\delta_1,b)$. Hence, $f(b)$ is the maximum value of $f$ on $[b-\delta_1, b]$. Let $\delta \in (0, \delta_1)$ be such that $|f(x)-f(b)| < \frac{\varepsilon}{2}$ for all $x \in [b-\delta, b]$, and let $f(x_1)$ take the minimum value of $f$ on $[b-\delta,b]$. Then $x_1<b$ and $f'_+(x_1)\ge 0$. Choose an arbitrary element $x_2\in (x_1,b)$, let $l(x) =f'_{-}(x_2)(x-x_2)+f(x_2)= cx+d$ with  $c=f'_{-}(x_2)$ and $d=-x_2f'_{-}(x_2)+f(x_2)$. Then $f(x) \geq l(x)$ for all $x \in [a,b]$. Since $c=f'_{-}(x_2)\geq f'_{+}(x_1) \geq 0$, we have
\begin{align*}
    f(x_0)&=f(b)=f(b)-f(x_2)+f(x_2) < \frac{\varepsilon}{2} + l(x_2) =\frac{\varepsilon}{2} + l(b)+l(x_2)-l(b) \\
         &=\frac{\varepsilon}{2}+l(x_0)+c(x_2-b)\le\frac{\varepsilon}{2}+l(x_0)<\varepsilon+l(x_0).
\end{align*}
This completes the proof.

\section{The generalization of Lemma~\ref{lem:convex equal}}\label{sec:generalization}

Observe that when $\mathfrak{A}=\mathbb{C}$, $\mathcal{S}(\mathfrak{A})$ contains only one point $\rho$ satisfying $\rho(z)=z$ for every $z\in\mathbb{C}$. In view of this observation,
a generalized version of Lemma~\ref{lem:convex equal} can be obtained in the Hilbert $C^*$-module case as follows.
\begin{theorem}\label{thm:conv func} Let  $t\in\mathcal{L}(E)_{\mbox{sa}}$ and $\rho\in \mathcal{S}(\mathfrak{A})$. Then for every continuous convex function  $f : [-\|t\|,\|t\|]\to \mathbb{R}$, we have
\begin{equation}\label{equ:inequality wrt convex function}f\big(\rho(\langle tx_0, x_0\rangle)\big)\leq \rho\big(\big\langle f(t)x_0, x_0\big\rangle\big),\end{equation}
where $x_0$ is an arbitrary element of $E$ satisfying $\rho(\langle x_0,x_0\rangle)=1$.
\end{theorem}
\begin{proof}
Given a positive linear functional $g$ on $\mathfrak{A}$, let $$\mathcal{N}_g=\{x\in E:g\left(\langle x,x\rangle\right)=0\}.$$
It follows from the Cauchy--Schwarz inequality that
$$\mathcal{N}_g=\{x\in E:g\left(\langle y,x\rangle\right)=g\left(\langle x,y\rangle\right)=0, \mbox{~for all~} y\in E\}.$$
Therefore, $\mathcal{N}_g$ is a closed linear subspace of $E$, and the quotient space $E/\mathcal{N}_g$ is a pre-Hilbert space equipped with the inner product $\langle\cdot,\cdot\rangle_g$ and the norm $\|\cdot\|_g$ defined by
\[
\langle x+\mathcal{N}_g,y+\mathcal{N}_g\rangle_g=g\big(\langle x,y\rangle\big),\quad \|x+\mathcal{N}_g\|_g=\sqrt{g\big(\langle x,x\rangle\big)}\quad (x,y\in E).
\]
Let $E_g$ be the completion of $E/\mathcal{N}_g$ and let $\iota_g:E\to E_g$ be the natural morphism, that is,
$$\iota_g(x)=x+\mathcal{N}_g\quad  (x\in E).$$

Suppose now that $t\in\mathcal{L}(E)_{\mbox{sa}}$ and $\rho\in \mathcal{S}(\mathfrak{A})$. For every $x,y\in E$, it is clear that
\begin{align*}
       |\rho(\langle tx,y\rangle)|^2
        & \leq \rho(\langle tx,tx\rangle)\rho(\langle y,y\rangle)=\rho(\langle t^2 x,x\rangle) \rho(\langle y,y\rangle) \\
        & \leq\rho(\langle \|t^2\|x,x\rangle)\rho(\langle y,y\rangle)=\|t\|^2\rho(\langle x,x\rangle)\rho(\langle y,y\rangle)\\&=\|t\|^2\cdot\|\iota_{\rho}(x)\|_{\rho}^2\cdot\|\iota_{\rho}(y)\|_{\rho}^2.
      \end{align*}
Hence
$$\varphi: E/\mathcal{N}_{\rho}\times E/\mathcal{N}_{\rho}\to\mathbb{C},\quad   \varphi\big(\iota_{\rho}(x),\iota_{\rho}(y)\big)=\rho(\langle tx,y\rangle)\quad (x,y\in E)$$
is a well-defined sesquilinear function  such that   for every $x,y\in E$,
$$\left|\varphi\big(\iota_{\rho}(x),\iota_{\rho}(y)\big)\right|\le \|t\|\cdot  \|\iota_{\rho}(x)\|_{\rho}\cdot\|\iota_{\rho}(y)\|_{\rho}$$
and
\begin{align*}\varphi\big(\iota_{\rho}(y),\iota_{\rho}(x)\big)&=\rho(\langle ty,x\rangle)=\rho(\langle y,tx\rangle)=\rho(\langle tx,y\rangle^*)=\overline{\varphi\big(\iota_{\rho}(x),\iota_{\rho}(y)\big)}.
\end{align*}
So, there exists a unique operator $T\in\mathbb{B}(E_{\rho})_{\mbox{sa}}$ such that
$\|T\|\le \|t\|$ and
$$\big\langle T\iota_{\rho}(x),\iota_{\rho}(y)\big\rangle_{\rho}=\rho(\langle tx,y\rangle)\quad (x,y\in E).$$
From the definition, we have
$$\big\langle \iota_{\rho}(tx),\iota_{\rho}(y)\big\rangle_{\rho}=\rho(\langle tx,y\rangle)\quad (x,y\in E).$$
This shows that for every $x,y\in E$,
$$\big\langle T\iota_{\rho}(x),\iota_{\rho}(y)\big\rangle_{\rho}=\big\langle \iota_{\rho}(tx),\iota_{\rho}(y)\big\rangle_{\rho}.$$
Therefore, by the density of $\mathcal{R}(\iota_{\rho})$ in $E_{\rho}$ it can be concluded that
$T\iota_{\rho}(x)=\iota_{\rho}(tx)$ for every $x\in E$, which in turn gives
$$T^n\iota_{\rho}(x)=\iota_{\rho}(t^nx),\quad \forall n\in\mathbb{N}\cup\{0\},  x\in E.$$
Consequently, for every polynomial $p$ and every $x\in E$ it is true that
\begin{equation}\label{equ:polynomial inequality}p(T)\big(\iota_{\rho}(x)\big)=\iota_{\rho}\big(p(t)x\big).\end{equation}
Obviously, $-\|t\|I_{E}\leq t\leq \|t\|I_{E}$. Also, due to $\|T\|\le \|t\|$ we have
$$-\|t\|I_{E_{\rho}}\leq T\leq \|t\|I_{E_{\rho}}. $$
So by \eqref{equ:polynomial inequality} and the continuity of $f$ on $[-\|t\|,\|t\|]$, it can be deduced that
\begin{equation}\label{equ:f inner product}f(T)\big(\iota_{\rho}(x)\big)=\iota_{\rho}\big(f(t)x\big)\quad (x\in E),\end{equation}
since there exists a sequence of polynomials $\{p_n\}_{n=1}^{\infty}$ such that $\{p_n\}_{n=1}^{\infty}$ converges to $f$ uniformly on $[-\|t\|,\|t\|]$.
By assumption we have
$$\|\iota_{\rho}(x_0)\|_{\rho}^2=\rho(\langle x_0,x_0\rangle)=1,$$
which leads by Lemma~\ref{lem:convex equal} to
\begin{equation}\label{equ:middle result01}f\big(\langle T\iota_{\rho}(x_0), \iota_{\rho}(x_0)\rangle_{\rho}\big)\leq \big\langle f(T)\iota_{\rho}(x_0),\iota_{\rho}(x_0) \big\rangle_{\rho}.\end{equation}
For every $x\in E$, by \eqref{equ:f inner product} we see that
\begin{align}\label{usefu101}&\big\langle T\iota_{\rho}(x), \iota_{\rho}(x)\big\rangle_{\rho}=\big\langle \iota_{\rho}(tx), \iota_{\rho}(x)\big\rangle_{\rho}=\rho(\langle tx,x\rangle),\\
\label{usefu102}&\big\langle f(T)\iota_{\rho}(x),\iota_{\rho}(x) \big\rangle_{\rho}=\big\langle \iota_{\rho}\big(f(t)x\big),\iota_{\rho}(x) \big\rangle_{\rho}=\rho\big(\big\langle f(t)x,x\big\rangle\big).
\end{align}
Hence, the desired inequality \eqref{equ:inequality wrt convex function} follows from \eqref{equ:middle result01}--\eqref{usefu102}.
\end{proof}

\begin{corollary}\label{cor1:conv func1} Let  $t\in\mathcal{L}(E)_+$ and $\rho\in \mathcal{S}(\mathfrak{A})$. Then \eqref{equ:inequality wrt convex function} is valid
for every continuous convex function  $f : [0,\|t\|]\to \mathbb{R}$ and every element $x_0$ of $E$ satisfying $\rho(\langle x_0,x_0\rangle)=1$.
\end{corollary}
\begin{proof}We follow the notations as in the proof of Theorem~\ref{thm:conv func}. By assumption $t\ge 0$, so it can be deduced from \eqref{usefu101} and \cite[Lemma~4.1]{Lance} that
$\big\langle T\iota_{\rho}(x), \iota_{\rho}(x)\big\rangle_{\rho}\ge 0$ for every $x\in E$. Since $\mathcal{R}(\iota_{\rho})$ is dense in $E_{\rho}$ and $T\in\mathbb{B}(E_{\rho})$, we have $T\in\mathbb{B}(E_{\rho})_+$. The desired conclusion can be derived by following the line in the proof of Theorem~\ref{thm:conv func}.
\end{proof}

\section{Some investigations on the H$\ddot{{\rm o}}$lder-McCarty inequalities}\label{sec:investigation}

In this section, we deal with the H$\ddot{{\rm o}}$lder-McCarty inequalities in the Hilbert $C^*$-module case. As a consequence of Corollary~\ref{cor1:conv func1}, we provide a known result\footnote{It is worth noting that \eqref{inequality for r bigger than 1} and
\eqref{inequality for r smaller than 1} were stated in \cite[Corollary~3.9]{MM} with $r-1$ being mistakenly written as $1-r$, and no proof of these two inequalities was provided therein.} as follows.
\begin{corollary}\label{cor:MM}{\rm \cite[Corollary~3.9]{MM}} Suppose that $t\in\mathcal{L}(E)_+$ and $\rho\in \mathcal{S}(\mathfrak{A})$. Then for every $x\in E\setminus\{0\}$,  we have
\begin{align}\label{inequality for r bigger than 1}&\big(\rho(\langle tx,x\rangle)\big)^r\leq\|x\|^{2(r-1)}\rho(\langle t^rx,x\rangle)\quad \mbox{for $r\geq 1$},\\
\label{inequality for r smaller than 1}&\big(\rho(\langle tx,x\rangle)\big)^r\geq\|x\|^{2(r-1)}\rho(\langle t^rx,x\rangle)\quad \mbox{for $0<r\le 1$}.
\end{align}
\end{corollary}
\begin{proof}We may as well assume that $r\ne 1$. If $\rho(\langle x,x\rangle)=0$, then by  the Cauchy--Schwarz inequality
$\rho(\langle y,x\rangle)=0$ for every $y\in E$, so in this case \eqref{inequality for r bigger than 1} and \eqref{inequality for r smaller than 1} are both trivially satisfied.
In what follows we assume that $\rho(\langle x,x\rangle)\ne 0$, which ensures that $\rho(\langle x_0,x_0\rangle)=1$, where $x_0$ is defined by
$$x_0=\frac{x}{\sqrt{\rho(\langle x,x\rangle)}}.$$

First, we consider the case that $r>1$. Let
$f(u)=u^r$ for $u\in [0,+\infty)$.
Since $r> 1$, $f$ is a continuous convex function on $[0,+\infty)$. So a direct application of Corollary~\ref{cor1:conv func1} yields
\begin{equation*}\big(\rho(\langle tx_0, x_0\rangle)\big)^r\leq \rho\big(\big\langle t^r x_0, x_0\big\rangle\big).\end{equation*}
It follows that
\begin{align*}\big(\rho(\langle tx, x\rangle)\big)^r\leq& \big(\rho(\langle x, x\rangle)\big)^{r-1}\rho\big(\big\langle t^r x, x\big\rangle\big)\le \|\langle x, x\rangle\|^{r-1}\rho\big(\big\langle t^r x, x\big\rangle\big)\\
=&\|x\|^{2(r-1)}\rho(\langle t^rx,x\rangle).\end{align*}
This shows the validity of \eqref{inequality for r bigger than 1}.

Next, we consider the case that $0<r<1$. In this case, the function $g$ defined by
$g(u)=-u^r$ is a continuous convex function on $[0,+\infty)$. Utilizing Corollary~\ref{cor1:conv func1} yields
\begin{equation*}\big(\rho(\langle tx_0,x_0\rangle)\big)^r   \geq \rho\big(\big\langle t^r x_0, x_0\big\rangle\big).\end{equation*}
Hence, \eqref{inequality for r smaller than 1} is derived as desired.
\end{proof}

We are now in the position to provide the Hilbert $C^*$-module version of the H$\ddot{{\rm o}}$lder-McCarty inequalities.
\begin{theorem}\label{thm:generalized McCarthy inequality} For every $t\in\mathcal{L}(E)_+$ and $x\in E\setminus\{0\}$,  we have
\begin{align}\label{inequality for r bigger than 3}&\|\langle tx,x\rangle\|^r\leq\|x\|^{2(r-1)}\|\langle t^rx,x\rangle\|\quad \mbox{for $r\geq 1$},\\
\label{inequality for r smaller than 4}&\|\langle tx,x\rangle\|^r\geq\|x\|^{2(r-1)}\|\langle t^rx,x\rangle\|\quad \mbox{for $0<r\le 1$}.
\end{align}
\end{theorem}
\begin{proof}For any normal element $a$ of $\mathfrak{A}$, by \cite[Theorem~3.3.6]{Murphy} there exists $\rho_0\in \mathcal{S}(\mathfrak{A})$ such that
$|\rho_0(a)|=\|a\|$. It follows that $\|a\|=\sup\limits_{\rho\in \mathcal{S}(\mathfrak{A})} |\rho(a)|$.

Now, suppose that $t\in\mathcal{L}(E)_+$ and $x\in E\setminus\{0\}$. In this case, by \cite[Lemma~4.1]{Lance} we see that
$\langle tx,x\rangle$ and $\langle t^rx,x\rangle$ are positive elements of $\mathfrak{A}$. So when $r\geq 1$,  we may use \eqref{inequality for r bigger than 1} to obtain
\begin{align*}\|\langle tx,x\rangle\|^r=&\sup\limits_{\rho\in\mathcal{S}(\mathfrak{A})}\big(\rho(\langle tx,x\rangle)\big)^r\leq\|x\|^{2(r-1)}\sup\limits_{\rho\in\mathcal{S}(\mathfrak{A})} \rho(\langle t^rx,x\rangle)\\
=&\|x\|^{2(r-1)}\|\langle t^rx,x\rangle\|.\end{align*}
This completes the proof of \eqref{inequality for r bigger than 3}. The proof of \eqref{inequality for r smaller than 4} is similar.
\end{proof}

It is known (see e.g.\,\cite[Proposition~1.3.5]{Pedersen}) that $\|a\|\le \|b\|$ whenever $a$ and $b$ are positive elements in a $C^*$-algebra satisfying $a\le b$.
This observation, together with Theorem~\ref{thm:generalized McCarthy inequality}, leads us  to investigate the following two issues.

\begin{problem}\label{prob01}Suppose that $\mathfrak{A}$ is a $C^*$-algebra. Is it true that
\begin{equation*}\langle tx,x\rangle^r\leq\|x\|^{2(r-1)}\langle t^rx,x\rangle
\end{equation*}
for every Hilbert $\mathfrak{A}$-module $E$, $t\in \mathcal{L}(E)_+$, $x\in E\setminus\{0\}$ and  $r>1$?
\end{problem}

\begin{problem}\label{prob02}Suppose that $\mathfrak{A}$ is a $C^*$-algebra. Is it  true that
\begin{equation*}\langle tx,x\rangle^r\geq\|x\|^{2(r-1)}\langle t^rx,x\rangle
\end{equation*}
for every Hilbert $\mathfrak{A}$-module $E$, $t\in \mathcal{L}(E)_+$, $x\in E\setminus\{0\}$ and  $0<r<1$?
\end{problem}

\begin{remark}\label{rem:Hilbert case}{\rm Every Hilbert module over the complex field is a Hilbert space, and thus the inequalities involved in Problems~\ref{prob01} and \ref{prob02} are reduced to the H$\ddot{{\rm o}}$lder-McCarty inequalities. So, the answers to
Problems~\ref{prob01} and \ref{prob02} are positive in the case that $\mathfrak{A}=\mathbb{C}$. Our next theorem shows that the same is true whenever $\mathfrak{A}$ is a commutative $C^*$-algebra.
}\end{remark}

\begin{theorem}\label{thm:commutative C-star alg case}For every commutative $C^*$-algebra $\mathfrak{A}$,  the answers to Problems~\ref{prob01} and \ref{prob02} are positive.
\end{theorem}
\begin{proof}For every self-adjoint element $a$ in a $C^*$-algebra $\mathfrak{A}$, it is known that $a\in \mathfrak{A}_+$ if and only if $\rho(a)\ge 0$ for every
$\rho\in  \mathcal{PS}(\mathfrak{A})$, where $\mathcal{PS}(\mathfrak{A})$ denotes the set consisting of all pure states on $\mathfrak{A}$.
When $\mathfrak{A}$ is commutative, by \cite[Theorem~5.1.6]{Murphy} we know that for any bounded linear functional $\varphi$ on $\mathfrak{A}$,
$\varphi\in \mathcal{PS}(\mathfrak{A})$ if and only if
$\varphi\ne 0$ and $\varphi$ is multiplicative (that is, $\varphi: \mathfrak{A}\to \mathbb{C}$ is a non-zero homomorphism), which ensures $\varphi(a^s)=\big(\varphi(a)\big)^s$ for every $a\in \mathfrak{A}_+$ and $s>0$, since there exists a sequence of polynomials $\{p_n(u)\}_{n=1}^\infty$ with $p_n(0)=0$ for each $n\in\mathbb{N}$ such that
$\{p_n(u)\}_{n=1}^\infty$ converges uniformly to the function $u^s$ on the interval $[0,\|a\|]$.

Suppose now that $\mathfrak{A}$ is a commutative $C^*$-algebra, $r>1$,  $E$ is a Hilbert $\mathfrak{A}$-module, $t\in \mathcal{L}(E)_+$ and $x\in E\setminus\{0\}$.
Let $a\in \mathfrak{A}_{\mbox{sa}}$ be given by
\begin{equation}\label{concrete self-adjoint element a}a=\|x\|^{2(r-1)}\langle t^rx,x\rangle-\langle tx,x\rangle^r.\end{equation}
For every $\rho\in \mathcal{PS}(\mathfrak{A})$, since $\rho$ is multiplicative we have
\begin{align*}\rho(a)=&\|x\|^{2(r-1)}\rho(\langle t^rx,x\rangle)-\rho(\langle tx,x\rangle^r)\\
=&\|x\|^{2(r-1)}\rho(\langle t^rx,x\rangle)-\big(\rho(\langle tx,x\rangle)\big)^r.
\end{align*}
Hence, by \eqref{inequality for r bigger than 1} we obtain $\rho(a)\ge 0$. The positivity of $a$ then follows from the arbitrariness of $\rho$ in $\mathcal{PS}(\mathfrak{A})$. This gives the positive answer to Problem~\ref{prob01}.
Similar reasoning yields the positive answer to Problem~\ref{prob02}.
\end{proof}

 Suppose that $\mathfrak{B}$  is a $C^*$-algebra. Let
$$\langle x,y\rangle=x^*y,\quad \mbox{for $x,y\in \mathfrak{B}$}.$$
With the inner-product given as above,
 $\mathfrak{B}$ is a Hilbert $\mathfrak{B}$-module,  and $\mathfrak{B}$ can be embedded into $\mathcal{L}(\mathfrak{B})$
via $b\to L_b$ \cite[Section~3]{LMX}, where
$L_b$ is defined by $L_b(x)=bx$ for $x\in \mathfrak{B}$. If  $\mathfrak{B}$ has a unit, then it can be easily derived from \eqref{equ:keep module operator} that $\mathcal{L}(\mathfrak{B})=\{L_b:b\in \mathfrak{B}\}$. In what follows, we identify $\mathfrak{B}$ with $\mathcal{L}(\mathfrak{B})$ for every unital $C^*$-algebra $\mathfrak{B}$.

With the convention described as above, we see that for every unital $C^*$-algebra $\mathfrak{A}$ considered as a Hilbert module over itself,  the inequalities involved in Problems~\ref{prob01} and \ref{prob02} turn out to be as follows:
 \begin{align}\label{csta alg r bigger than 1}&(x^*tx)^r\leq\|x\|^{2(r-1)}x^* t^r x,\quad \forall t\in \mathfrak{A}_+, x\in \mathfrak{A}, r>1,\\
 \label{csta alg r smaller than 1}&(x^*tx)^r\geq\|x\|^{2(r-1)}x^* t^r x,\quad \forall t\in \mathfrak{A}_+, x\in \mathfrak{A}, 0<r<1.
\end{align}

 \begin{lemma}\label{lem:counterexample01} Let $\mathfrak{A}=M_2(\mathbb{C})$. Then \eqref{csta alg r bigger than 1} fails to be true for
 \begin{equation*}x=\left(
                      \begin{array}{cc}
                        1 & 1 \\
                        0 & 1 \\
                      \end{array}
                    \right).\quad t=\left(
                                      \begin{array}{cc}
                                        2 & 1 \\
                                        1 & 2 \\
                                      \end{array}
                                    \right),\quad r=3.
 \end{equation*}
 \end{lemma}
 \begin{proof}For the given $x,t$ and $r$,  let $A$ and $B$ be the left side and the right side of \eqref{csta alg r bigger than 1}, respectively. A direct use of the Matlab shows that
 \begin{align*}&A=(x^*tx)^3=\left(
                              \begin{array}{cc}
                                98 & 183\\
                                183 & 342 \\
                              \end{array}
                            \right),\\
 &B=\|x\|^{-4/3}x^* t^{1/3} x\approx \left(
                                     \begin{array}{cc}
                                       95.9574  & 185.0608 \\
                                       185.0608 & 370.1215\\
                                     \end{array}
                                   \right),\\
 &C:=B-A\approx \left(
                           \begin{array}{cc}
                             -2.0426  &   2.0608\\
                              2.0608 &  28.1215 \\
                           \end{array}
                         \right).
 \end{align*}
 Since $-2.0426<0$, we see that the matrix  $C$ is not semi-positive definite.
 \end{proof}

 \begin{lemma}\label{lem:counterexample02} Let $\mathfrak{A}=M_2(\mathbb{C})$. Then \eqref{csta alg r smaller than 1} fails to be true for
 \begin{equation*}x=\left(
                      \begin{array}{cc}
                        9 & 9 \\
                        1 & -25 \\
                      \end{array}
                    \right)
 ,\quad t=\left(
             \begin{array}{cc}
               125 & 75 \\
               75 & 45 \\
             \end{array}
           \right)
 ,\quad r=1/4.\end{equation*}
 \end{lemma}
 \begin{proof}For the given $x,t$ and $r$,  let $A$ and $B$ be the left side and the right side of \eqref{csta alg r smaller than 1}, respectively. Direct computations yield
 \begin{align*}&A=(x^*tx)^{1/4}\approx \left(
                              \begin{array}{cc}
                                 8.0901 & -5.0563\\
                                -5.0563 & 3.1602\\
                              \end{array}
                            \right),\\
 &B=\|x\|^{-3/2}x^* t^{1/4} x\approx \left(
                                     \begin{array}{cc}
                                       1.7772  &  -1.1105 \\
                                        -1.1105 & 0.6956\\
                                     \end{array}
                                   \right),\\
 &C:=A-B\approx \left(
                           \begin{array}{cc}
                             6.3130 &   -3.9458\\
                              -3.9458 &   2.4646 \\
                           \end{array}
                         \right).
 \end{align*}
 It follows that $\mbox{det}(C)\approx -0.0108<0$, where $\mbox{det}(C)$ denotes the determinant of $C$. So, the matrix $C$ is not semi-positive definite.
 \end{proof}

A combination of Lemmas~\ref{lem:counterexample01} and \ref{lem:counterexample02} yields the following theorem.
\begin{theorem}\label{thm:negative answer}If $H$ is a Hilbert space with $\mbox{dim}(H)\ge 2$, then the answers to Problems~\ref{prob01} and \ref{prob02} are negative for the $C^*$-algebras $\mathbb{B}(H)$ and $\mathbb{K}(H)$, where $\mathbb{K}(H)$ denotes the set consisting of all compact operators on $H$.
\end{theorem}

\vspace{2ex}






\begin{thebibliography}{99}

\bibitem{BDMP}P. Bhunia, S. S. Dragomir, M. S. Moslehian and K. Paul, Lectures on numerical radius inequalities, Infosys Science Foundation Series in Mathematical Sciences, Springer, Cham, 2022.





\bibitem{KR}F. Kittaneh and  M. H. M. Rashid, Inequalities for norms and numerical radii of operator
matrices, Rev. Real Acad. Cienc. Exactas Fis. Nat. Ser. A-Mat. (2024) 118:150.

\bibitem{Lance}E. C. Lance, Hilbert $C^*$-modules--A toolkit for operator algebraists, Cambridge University Press, Cambridge,  1995.

\bibitem{LMX}W. Luo, M. S. Moslehian and Q. Xu, Halmos' two projections theorem for Hilbert $C^*$-module operators and the Friedrichs angle of two closed submodules, Linear Algebra Appl. 577 (2019), 134--158.




\bibitem{MT} V. M. Manuilov and E. V. Troitsky, Hilbert $C^*$-modules, Translated from the 2001 Russian original by the authors, Translations of Mathematical Monographs, 226, American Mathematical Society, Providence, RI, 2005.

\bibitem{McCarthy}C. A. McCarthy,  $C_p$,  Israel. J. Math. 5 (1967), 249--271.

\bibitem{MRO}A. K. Mirmostafaee, O. P.  Rahpeyma and M. E. Omidvar, Numerical radius inequalities for finite sums of operators, Demonstr. Math. 47 (2014), no. 4, 963--970.

\bibitem{MM}S. F. Moghaddam and A. K. Mirmostafaee, Numerical radius inequalities for Hilbert $C^*$-modules, Math. Bohem. 147 (2022), no. 4, 547--566.




\bibitem{MP}B. Mond and J. E. Pe$\breve{{\rm c}}$ari$\acute{{\rm c}}$,  Convex inequalities in Hilbert space, Houston J. Math. 19 (1993), no. 3, 405--420.





\bibitem{Murphy} G. J. Murphy, $C^*$-algebras and operator theory, Academic Press, New York, 1990.

\bibitem{Paschke}W. L. Paschke, Inner product modules over $B^*$-algebras, Trans. Amer. Math. Soc. 182 (1973), 443--468.


\bibitem{Pedersen} G. K. Pedersen, $C^*$-algebras and their automorphism groups (London Math. Soc. Monographs 14), Academic Press, New York, 1979.


\bibitem{Rashid}M. H. M. Rashid, Power inequalities for the numerical radius of
operators in Hilbert spaces, Khayyam J. Math. 5 (2019), no. 2, 15--29.

\bibitem{RA}M. H. M. Rashid and N. H. Altaweel, Some generalized numerical radius
inequalities for Hilbert space operators, J. Math. Inequal. 16 (2022), no. 2, 541--560.

\bibitem{Yang-Li}C. Yang and D. Li, Some improvements about numerical radius
inequalities for Hilbert space operators, J. Math. Inequal. 18 (2024), no. 1, 219--234.









\end{thebibliography}
\end{document}